\documentclass[a4paper,10pt]{article}
\usepackage{amsmath,amsfonts,amssymb,amstext,amsthm,amscd}
\usepackage{amsmath,amssymb,latexsym,amscd,fancybox,graphicx}
\usepackage[all]{xy}
\newtheorem{prop}{Proposition}
\newtheorem{theo}{Theorem}
\newtheorem{defi}{Definition}

\title{Pfaffian Systems from Twistor Fibrations\footnote{
Pfaffian system, twistor fibration, superhorizontal distribution,
infinitesimal symmetries.} }
\author{Ramiro Carrillo-Catal\'{a}n \thanks{The first author is supported by a Japanese Government (Monbukagakusho: MEXT) Scholarship.} \and {Akio Yabe }}
\date{September 11, 2006}
\begin{document}
\maketitle

\begin{abstract}
Canonical twistor fibrations lead to Pfaffian systems by means of
their superhorizontal distribution. The aim of this note is to
identify explicitly the Pfaffian systems of five or less variables
that arise in this way in terms of  the classification given in
[1].
\end{abstract}

 \section{Introduction}

\ A Pfaffian system on an $n$-dimensional differentiable manifold
is a collection of
 1-forms which are linearly independent at each point. Geometrically this system gives
 rise  to a subbundle of the tangent space, that is a
  distribution. Conversely a distribution induces an equivalence class of Pfaffian
 systems (see [1],[3],[6] or [9]).  For $n\le 5$ such systems were investigated in detail by E.
 Cartan.\\

 On the total space of any canonical twistor fibration (see [2] or [5]), there is a natural
 \lq\lq superhorizontal\rq\rq \ distribution, which is important in the theory of harmonic
 maps. The total space is a  flag manifold, and a corresponding Pfaffian system on the big cell of the
 flag manifold  can be constructed by using Lie-theoretic local
 coordinates to express integral curves of the distribution. The equations for these holomorphic curves are solved explicitly in [4].\\

In this paper we identify explicitly all  Pfaffian systems with
$n\le 5$ which arise from twistor fibrations in this way, in terms
of the known classification. It is of interest to consider the
significance of these twistor Pfaffian systems, and we discuss
briefly one aspect here, namely the Lie algebra of infinitesimal
symmetries.  Cartan observed that the Lie algebra of the
exceptional group $G_2$ can be characterized as the symmetry
algebra of such a system.  However, all the other systems with
$n\le 5$  arising from twistor fibrations have
infinite-dimensional symmetry algebras. In the
simplest example, we present an explanation for this phenomenon.\\

The paper  is organized as follows.  In section 2 we describe the
basic definitions and theorems concerning Pfaffian systems. In
addition, the classification of low dimensional Pfaffian systems
from [1] is reproduced at the end of this section (tables 1-3).
Section 3  reviews some facts from structure theory of Lie
algebras as well as the concept of twistor fibration.   The
construction of Pfaffian systems from twistor fibrations is
explained in section 4. Subsections 4.1 - 4.3 contain detailed
computations for three representative examples:
$\pi:F_3(\mathbb{C}) \rightarrow \mathbb{C}P^2$,
$\pi:SO(5)/U(1)\times SO(2) \rightarrow SO(5)/SO(2)\times SO(3)$
and $\pi:G_2/U(2) \rightarrow G_2/SO(4)$. The complete list of
results is given in table 4. Finally in the last section we
discuss the Lie algebra of infinitesimal symmetries.

\section{Pfaffian systems}

In this section we will review some basic concepts and facts on
Pfaffian systems theory (see [3]). Let $M$ be an $n-$dimensional
differentiable manifold. Let $TM$ and $T^{*}M$ denote the tangent
and cotangent bundle respectively.\\

  In the following $
\Omega^{k}(M)$ will denote the space of all differentiable k forms
on $M$. We have $
\Omega^{k}(M)=\Gamma(\bigwedge^{k}T^{*}M)$, where
$\Gamma(\bigwedge^{k}T^{*}M)$ denotes the space of sections of
$\bigwedge^{k}T^{*}M$.\\

A (local) \textit{Pfaffian system} of rank $r$ is a set:
\begin{align*}
S=\{\omega_{1}\dots\omega_{r}\} \qquad  \omega_{i}\in
\Omega^{1}(M)|_{U}.
\end{align*}
where the r differential 1-forms are defined on an open subset $U$
of $M$ and are required to be linearly independent at each point.\\

  Given a Pfaffian system $S$, the  vector subbundle $D$
  of $TM|_{U}$ with $n-r$ dimensional fibers $D_{p}=\bigcap_{i=1}^{r}\ker\omega_{i}(p)$ for all $\omega_{i}\in S$, $ p\in U\subset
 M$ is a \textit{distribution} on $U$. Conversely, a distribution
 give rise to an equivalence class of Pfaffian systems, since the annihilator subbundle of $D$, with fibers
\begin{align*}
D^{\perp}_{p}&=\{\omega\in T^{*}_{p}M | \  \omega(X)=0 \ \forall
X\in D_{p}\},
\end{align*}
is locally spanned by $r$ linearly independent $1$-forms
\textit{i.e.}\
$D_{p}^{\perp}=\text{span}\{(\omega_{1})_{p}\dots(\omega_{r})_{p}\}\subset
T^{*}_{p}M$.  The ideal of $\Gamma(\bigwedge^{*}T^{*}M)$ generated
by the system is:
\begin{align*}
\mathcal{I}(D)=
\bigoplus^{n}_{k=1}\mathcal{I}^{k}(D)=\Big\{\sum_{i=1}^{r}\alpha_{i}\wedge
\omega_{i} \ | \ \omega_{i}\in S, \ \alpha_{i}\in
\Omega^{*}(M)\Big\}.
\end{align*}
with  $\mathcal{I}^{k}(D)=\{\omega\in\Omega^{k}(M) \ | \
\omega(X_{1}\dots X_{k})=0 \quad \forall X_{i}\in D\}$ for $k\geq
1$.\\

At a point $p\in M$ the \textit{characteristic space} of
$\mathcal{I}(D)$ is defined as:
\begin{align*}
\mathcal{A}(\mathcal{I}(D))_{p}&=\{X_{p}\in T_{p}M  \ | \
X_{p}\rfloor \mathcal{I}(D)_{p}\subset \mathcal{I}(D)_{p} \}\\
&=\{ X_{p}\in T_{p}M  \ | \ \omega_{i}(X_{p})=0, \ \
X_{p}\rfloor\text{d}\omega_{i}\equiv 0 \mod S,  \ i=1\dots r \}.
\end{align*}
Its annihilator
$C(\mathcal{I}(D))_{p}:=\mathcal{A}(\mathcal{I}(D))_{p}^{\perp}\subset
T^{*}_{p}M$,  called the \textit{retracting } subspace at $p$, is
the smallest subspace of $T^{*}_{p}M$ such that
$\bigwedge^{*}(C(\mathcal{I}(D))$ contains a finite set $S$ of
elements generating $\mathcal{I}(D)$ as an ideal.
 The dimension of $ C(\mathcal{I}(D))_{p}$ is called the
\textit{class} of $S$ at $p$.\\

A \textit{symmetry} of $S$ is a local biholomorphism $\phi$ on $U$
which preserves $D$ (or $S$), \textit{i.e.}\ $\phi_{*}(D)=D$ (or
$\phi^{*}(S)=S$). An \textit{infinitesimal  symmetry} is a vector
field $X$ over $U$ which generates a local one parameter
transformation group of symmetries of $S$.  Namely $X$ satisfies:
$[X, \Gamma(D)]\subseteq \Gamma(D)$, or equivalently:
$L_{X}\omega_{i}\equiv 0 \ \mod S $ with $\omega_{i}\in S$. The
set of all infinitesimal symmetries of $S$ is denoted
 $L(S)$.\\

Then the characteristic space
$\mathcal{A}(\mathcal{I}(D))=\bigcup_{p\in U}
\mathcal{A}(\mathcal{I}(D))_{p}$ of $\mathcal{I}(D)$ can be
described (\textit{cf}.\ [6]) in terms of infinitesimal symmetries
as:
\begin{align}
\mathcal{A}(\mathcal{I}(D))=L(S)\cap\Gamma(D).
\end{align}

A submanifold $N$ of $M$ is called an \textit{integral manifold}
of $S$ (or $D$) if $T_{p}N \subseteq D_{p}$ for any $p\in N$.
Equivalently an integrable manifold is given by an immersion
$i:N\longmapsto M$
such that  $i^{*}(\alpha)=0$ for all $\alpha\in \mathcal{I}(D)$.\\

A Pfaffian system $S$ (or the distribution $D$) is said to be
\textit{completely integrable} if there exists an
integrable manifold of dimension $dim(M)-rank(S)$, passing through every point.\\

The integrability of Pfaffian systems can be characterized in
various ways:

\begin{theo}[Frobenius] Let $M$ be a differentiable manifold. Let $S=\{\omega_{i}\}^{r}_{i=1}$ by a
Pfaffian system of rank r with induced $n-r$ dimensional
distribution $D$ defined over $U\subset M$. Let  $\mathcal{I}(D)$
 be the ideal generated by the system as above. Then the following
are equivalent:
\begin{enumerate}
\item[i)] $S$ (or $D$) is completely integrable.
 \item[ii)] $D$ is involutive: \ $[\Gamma(D),\Gamma(D)]\subset
 \Gamma(D)$.
\item[iii)]$\mathcal{I}(D)$ is a differential ideal:
$\text{d}\mathcal{I}(D)\subset\mathcal{I}(D)$, \textit{i.e.}\
 $d\omega_{i}=\sum_{j=1}^{r}\alpha^{i}_{j}\wedge\omega_{j}$ for $\omega_{j}\in S$ and  $\alpha^{i}_{j}\in\Omega^{k}(M)|_{U}$.
  Briefly  we write:  $d\omega_{i}\equiv 0 \mod S$.
  \item[iv)] Locally there exists
a coordinate system $(y_{1}\dots y_{n})$ such that
$\mathcal{I}(D)$ is generated
 by $(dy_{1}\dots dy_{r})$.
 \item[v)] $\text{rank}_{p}(S)=\text{class}_{p}(S)$ for every $p\in U\subset
 M$. In this case
$\mathcal{A}(\mathcal{I}(D))=D$ and $C(\mathcal{I}(D))=S$.\\
\end{enumerate}
\end{theo}

 A natural problem in  Pfaffian systems theory is the
classification  of systems with constant rank, which requires the
study of invariants. The class and the rank are
 fundamental invariants that entirely determine the local
classification of completely integrable systems by the Frobenius theorem.\\

For a nonintegrable Pfaffian system $S$ of rank $r$ and
corresponding distribution $D$, we define the subbundle $D'$ which
has  fibers spanned by all elements of $T_{p}M$ of the form
$X_{p}+[Y_{p},Z_{p}] $ with $X,Y,Z\in\Gamma(D)$, \textit{i.e.}\
\begin{align*}
D'=\Gamma(D)+[\Gamma(D), \Gamma(D)].
\end{align*}
 Notice that $D_{p}\subset D'_{p}$. A representative $S'$ of the equivalence class of Pfaffian systems corresponding to this distribution is called (with abuse of notation) the \textit{derived system} of
 $S$.  This leads to the notion of \textit{higher derived systems}. The $k$-th derived
system $S^{k}$ is defined successively by:
\begin{align*}
S^{k+1}=(S^{k})',  \ \  \  \ S^{0}:=S .
\end{align*}

It follows that the rank and the class of $S^{k}$ are invariants
of such a Pfaffian systems.   $S$ is called \textit{regular} if
all derived systems are of constant rank. For a regular system
$S$, there exists $\mu\in \mathbb{Z}_{\ge 0}$ such that:
\begin{align*}
S^{\mu+1}=S^{\mu}\subset\dots \subset S^{k}\subset S^{k-1}\subset
\dots \subset S''\subset S'\subset S.
\end{align*}

$S^{\mu}$ is always 0 or completely integrable; It is the
smallest completely integrable subsystem of $S$. (See [10],[1]).\\

The basic idea to achieve the classification of nonintegrable
Pfaffian systems of constant rank is to look for completely
integrable subsystems of S and use the defined
invariants to determine a local model of the system.\\

In [1] the equivalence classes  of low dimensional ($n \leq 5$)
Pfaffian systems of constant rank and class is presented by giving
 local models of those systems. Nevertheless this is a partial classification, since for the cases
$S^{11}_{5}(f)$ and $S^{15}_{5}(f)$ the local equivalence classes
are not specified. For the sake of reference this classification
is reproduced in the following tables.
\begin{itemize}
\item [\textbf{Table 1.}] Pfaffian systems of dimension 3.

\begin{center}
    \begin{tabular}{|c|c|c|}            \hline
        Pfaffian system $ S $   & $\mathrm{rank}(S)$    & $\mathrm{class}(S)$   \\ \hline
        $ S^1_3 = \left\{ dz_1 \right\} $   & 1 & 1 \\ \hline
        $ S^2_3 = \left\{ dz_1+z_2dz_3 \right\} $ & 1 & 3 \\ \hline
        $ S^3_3 = \left\{ dz_1 ,~ dz_2 \right\} $ & 2 & 2 \\ \hline
        $ S^4_3 = \left\{ dz_1 ,~ dz_2 ,~ dz_3 \right\} $ & 3 & 3 \\ \hline
    \end{tabular}
\end{center}

\item [\textbf{Table 2.}]  Pfaffian systems of dimension 4.

\begin{center}
    \begin{tabular}{|c|c|c|c|c|}            \hline
        Pfaffian system $ S $   & $\mathrm{rank}(S)$    & $\mathrm{class}(S)$   & $\mathrm{rank}(S^\prime)$ & $\mathrm{class}(S^\prime)$    \\ \hline
        $ S^1_4 = \left\{ dz_1 \right\} $   & 1 & 1 &  &  \\ \hline
        $ S^2_4 = \left\{ dz_1 ,~ dz_2 ,~ dz_3 \right\} $ & 3 & 3 &  & \\ \hline
        $ S^3_4 = \left\{ dz_1 ,~ dz_2 ,~ dz_3 ,~ dz_4  \right\} $ & 4 & 4 &  &  \\ \hline
        $ S^4_4 = \left\{ dz_1+z_2dz_3 \right\} $ & 1 & 3 &  &  \\ \hline
        $ S^5_4 = \left\{ dz_1 ,~ dz_2 \right\} $ & 2 & 2 &  &  \\ \hline
        $ S^6_4 = \left\{ dz_1 ,~ dz_2+z_3dz_4 \right\} $ & 2 & 4 & 1 & 1 \\ \hline
        $ S^7_4 = \left\{ dz_1+z_2dz_3 ,~ dz_2+z_4dz_3 \right\} $ & 2 & 4 & 1 & 3 \\ \hline
    \end{tabular}
\end{center}
\newpage
\item [\textbf{Table 3.}]   Pfaffian systems of dimension 5 .
\end{itemize}
\begin{center}
    \rotatebox{270}{$
    \begin{tabular}{|c|c|c|c|c|c|c|}            \hline
        Pfaffian~system~S   & $\mathrm{rank}(S)$    & $\mathrm{class}(S)$   & $\mathrm{rank}(S^\prime)$ & $\mathrm{class}(S^\prime)$    & $\mathrm{rank}(S^{\prime\prime})$ & $\mathrm{class}(S^{\prime\prime})$    \\ \hline
        $ S^1_5 = \left\{ dz_1 \right\} $   & 1 & 1 &&&&\\ \hline
        $ S^2_5 = \left\{ dz_1+z_2dz_3 \right\} $ & 1 & 3 &&&&\\ \hline
        $ S^3_5 = \left\{ dz_1+z_2dz_3+z_4dz_5 \right\} $ & 1 & 5 &&&&\\ \hline
        $ S^4_5 = \left\{ dz_1 ,~ dz_2 ,~ dz_3 ,~ dz_4 , dz_5 \right\} $ & 5 & 5 &&&&\\ \hline
        $ S^5_5 = \left\{ dz_1 ,~ dz_2 ,~ dz_3 ,~ dz_4 \right\} $ & 4 & 4 &&&&\\ \hline
        $ S^6_5 = \left\{ dz_1 ,~ dz_2 ,~ dz_3 \right\} $ & 3 & 3 &&&&\\ \hline
        $ S^7_5 = \left\{ dz_1 ,~ dz_2 ,~ dz_3+z_4dz_5 \right\} $ & 3 & 5 &2&2&&\\ \hline
        $ S^8_5 = \left\{ dz_1 ,~ dz_2+z_3dz_4 ,~ dz_3+z_5dz_4 \right\} $ & 3 & 5 &2&4&1&1\\ \hline
        $ S^9_5 = \left\{ dz_1+z_2dz_3 ,~ dz_2+z_4dz_3 ,~ dz_3+z_5dz_4 \right\} $ & 3 & 5 &2&4&1&3\\ \hline
        $ S^{10}_5 = \left\{ dz_1+z_2dz_3 ,~ dz_2+z_4dz_3 ,~ dz_4+z_5dz_3 \right\} $ & 3 & 5 &2&4&1&3\\ \hline
        $ S^{11}_5(f) = \left\{ dz_1+z_3dz_4 ,~ dz_2+z_5dz_3+fdz_4 ,~ dz_3+\frac{\partial{f}}{\partial{z_5}}dz_4 \right\} $ & 3 & 5 &2&5&&\\ \hline
        $ S^{12}_5 = \left\{ dz_1 ,~ dz_2 \right\} $ & 2 & 2 &&&&\\ \hline
        $ S^{13}_5 = \left\{ dz_1,~ dz_2+z_3dz_4 \right\} $ & 2 & 4 &1&1&&\\ \hline
        $ S^{14}_5 = \left\{ dz_1+z_2dz_3 ,~ dz_2+z_4dz_3 \right\} $ & 2 & 4 &1&3&&\\ \hline
        $ S^{15}_5(f) = \left\{ dz_1+z_3dz_4 ,~ dz_2+z_5dz_3+fdz_4 \right\} $ & 2 & 5 &&&&\\ \hline
    \end{tabular}$}
    \end{center}


\section{Twistor fibrations}

  The aim of this section is to review the twistor fibration concept.  We
  start by  recalling some facts from the  structure theory of Lie algebras.\\

Let  $\mathfrak{g}$ be a Lie algebra of a compact semisimple Lie
group with Cartan subalgebra $\mathfrak{t}$. Let
$\mathfrak{g}^\mathbb{C}$
be the complexification of  $\mathfrak{g}$, and  $\mathfrak{t}^\mathbb{C}$ that of  $\mathfrak{t}$.\\

 A functional
$\alpha\in(\mathfrak{t}^{\mathbb{C}})^{*}$ is called \textit{a
root} of $\mathfrak{g}^{\mathbb{C}}$ (with respect to
$\mathfrak{t}^{\mathbb{C}}$) if $\mathfrak{g}_{\alpha}\neq \{0\},
\ $where $ \
\mathfrak{g}_{\alpha}=\{X\in\mathfrak{g}^{\mathbb{C}}\ | \
\text{ad}(H)X=\alpha(H)X\ \text{for all} \
H\in\mathfrak{t}^{\mathbb{C}} \} $ is the
 \textit{ root space } of $\alpha$.  The  Lie algebra $\mathfrak{g}^{\mathbb{C}}$ has then a root decomposition
into the Cartan subalgebra and the root spaces:
\begin{align*}
\mathfrak{g}^{\mathbb{C}}=\mathfrak{t}^{\mathbb{C}}\oplus\big(\bigoplus_{\alpha\in\Delta}\mathfrak{g}_{\alpha}\big).
\end{align*}
Here $\Delta$ denotes the set of roots of
$\mathfrak{g}^{\mathbb{C}}$. This is a finite subset of
$(\mathfrak{t}^{\mathbb{C}})^{*}$.\\

It is also necessary to recall the existence of  a subset
$\Pi=\{\alpha_{1}\dots\alpha_{r}\}$
($r=\text{dim}_{\mathbb{C}}\mathfrak{t}^{\mathbb{C}}$) of $\Delta$
such that every root can be expressed uniquely as
$\alpha=\sum_{i=1}^{l}n_{i}\alpha_{i}$ where the $n_{i}$ are
integers, either all nonnegative ($\alpha$ is then a
\textit{positive} root) or all non-positive ($\alpha$ is a
\textit{negative} root).  Such a set is called \textit{set of
simple roots} for $\Delta$. Each
 subset $I=\{\alpha_{i_{1}},\dots,\alpha_{i_{s}}\}$ of $\Pi$ defines on $\Delta$
a \textit{height} function $n_{I}$ by $n_{I}(\alpha)=\sum_{i\in I}n_{i}$.\\

According to this any parabolic subalgebra  of
$\mathfrak{g}^{\mathbb{C}}$ can be expressed as:
\begin{align*}
\mathfrak{p}_{I}=\mathfrak{t}^{\mathbb{C}}\oplus\Big(\bigoplus_{i
\leq 0}\mathfrak{g}_i^I
\Big),~~~\text{with}~~~\mathfrak{g}_i^I=\bigoplus_{n_I(\alpha)=i}\mathfrak{g}_{\alpha}.
\end{align*}

Given a parabolic subalgebra  $\mathfrak{p}_{I}$, there exists a
unique element $\xi\in \mathfrak{t}^{\mathbb{C}}$ such that
$\alpha(\xi)=\sqrt{-1}n_{I}(\alpha)$. \textit{i.e.}\
$$\bigoplus_{\alpha(\xi)=\sqrt{-1}i}
\mathfrak{g}_{\alpha}=\mathfrak{g}_i^{I}.$$
 This element is
called the \textit{canonical element} of $\mathfrak{p}_{I}$.  This
 element  defines an involutive automorphism $\tau_{\xi}=\text{Ad}(\exp(\pi\xi))$ of $\mathfrak{g}$ (the complex linear extension to $\mathfrak{g}^{\mathbb{C}}$)  called the
\textit{canonical involution}  for $\mathfrak{p}_{I}.$\\

 The facts mentioned above are intimately related with the concept
 of \textit{generalized flag manifolds}, which are homogeneous
 spaces of the form $G/H$ where $G$ is a compact Lie group and $H$ is the
centralizer of a torus in $G$. We have a natural isomorphism
 $G^{\mathbb{C}}/P_{I}\cong G/H $ where $G^{\mathbb{C}}$ is the complexification of $G$
and $P_{I}= \{g\in G^{\mathbb{C}}|
\text{Ad}(g)\mathfrak{p}_{I}\subseteq \mathfrak{p}_{I}\}$ is a
parabolic subgroup of $G^{\mathbb{C}}$ .\\

 Let  $G/H\cong
G^{\mathbb{C}}/P_{I}$ be a generalized flag manifold and let $\xi$
be the canonical element of the Lie algebra $\mathfrak{p}_{I}$ of
$P_{I}$ with canonical automorphism $\tau_{\xi}$.  Then the
following is true:

\begin{itemize}
\item $G^{\mathbb{C}}/P_{I}\cong \text{Ad}(G)\xi$.

\item $H=P_{I}\cap G$.

 \item For $g=\exp\pi\xi$ define $\sigma_{g}=R_{g^{-1}}\circ L_{g}$.  If $K$ is a closed  subgroup of
  $G^{\mathbb{C}}$ such that  $H\subseteq K\subseteq G$, \
 $G_{\sigma_{g}}^{0}\subseteq K\subseteq G_{\sigma_{g}}$, then  $G/K$ is a \textit{symmetric
 space}. Here $ G_{\sigma_{g}}$ denotes the set of fixed points of $\sigma_{g}$ and $G_{\sigma_{g}}^{0} $ its
 identity component.
\end{itemize}

Finally we have:\\

\begin{defi}
The natural map $\pi_{\xi}:G/H \longmapsto G/K$ is called the
canonical \emph{(Burstall-Rawnsley})  \textit{twistor fibration}
associated to the generalized flag manifold $G/H\cong
G^{\mathbb{C}}/P_{I}$.\\
\end{defi}

For details we refer to  [5], [2] and [7].\\

\section{Pfaffian systems from twistor fibrations}

\   Now we are in a position to describe how Pfaffian systems and
twistor fibrations are related.  The key is to consider a certain
subbundle of the tangent bundle of the
generalized flag manifold of the twistor fibration.\\

Let $G^{\mathbb{C}}/P_{I}$ be a generalized flag manifold. Recall
that as an homogeneous space its  tangent bundle is isomorphic to
$G^{\mathbb{C}}\times_{P_I}\left(\mathfrak{g}^{\mathbb{C}}/\mathfrak{p}_I\right)$.
An explicit isomorphism
$f:G^{\mathbb{C}}\times_{P_I}\left(\mathfrak{g}^{\mathbb{C}}/\mathfrak{p}_I\right)\rightarrow
TG^{\mathbb{C}}/P_I$ is given by:
\begin{align}
f([(g,[Y])])=\left.\frac{d}{ds}g \circ \exp sY  P_I \right|_{s=0}.
\label{isomorfismo}
\end{align}

\begin{defi}
The subbundle $SH$ of $TG^{\mathbb{C}}/P_I$ which corresponds to
\begin{align*}
G^{\mathbb{C}}\times_{P_I}\left((\mathfrak{p}_I\oplus\mathfrak{g}^{I}_{1})/\mathfrak{p}_I\right)
\end{align*}
is called the \textbf{superhorizontal distribution}.
\end{defi}

This concept arises from  the study of harmonic maps  that come from twistor fibrations. \\

We shall construct a particular Pfaffian system corresponding to
the superhorizontal distribution by expressing integral curves of
the distribution in Lie theoretic local coordinates.  The equations for integral curves of this distribution are solved expliclitly in [4], (see page 561).

 It is well known  that local coordinates for the generalized flag manifold, $G^{\mathbb{C}}/P_{I}$ are provided by the
 ``\textit{big cell}'' $[\exp \bigoplus_{i>0}\mathfrak{g}^{I}_{i}]\cong \mathbb{C}^{m}$.  Locally therefore we can construct any complex integral curve $\Phi:\mathbb{C}\longmapsto G^{\mathbb{C}}/P_{I}$  by
 giving a holomorphic map
 $C:\mathbb{C}\longmapsto\bigoplus_{i>0}\mathfrak{g}_{i}^{I}$ and setting $\Phi=\rho(\exp C)$, where $\rho$ denotes
 the projection  $G^{\mathbb{C}}\longmapsto G^{\mathbb{C}}/P_{I}$.\\

  The condition for this integral curve to be tangent to the superhorizontal distribution, \textit{i.e.}\ $\text{d}\Phi \big
(\frac{\partial}{\partial z}\big )\in SH$, is that:
\begin{align}
\left((\exp C)^{-1}(\exp C)^\prime \right)_{\mathfrak{g}^{I}_i}= 0
~~~~~~(i \geq 2)\label{condition}
\end{align}
In other words, $(\exp -C)(\exp C)^\prime$  should have  zero
component in $\bigoplus_{i \geq 2}\mathfrak{g}^{I}_i$.\\

 Using the formula for the derivative of the exponential map one has:
\begin{align*}
(\exp -C)(\exp C)^\prime &= \frac{I-e^{-\text{ad}C}}{\text{ad}C} C^\prime \\
            &=C^\prime-\frac{1}{2!}\text{ad}C C^\prime +\frac{1}{3!}(\text{ad}C)^2C^\prime-\cdots \\
            &=C^\prime-\frac{1}{2}[C,C^\prime]+\frac{1}{6}[C,[C,C^\prime]]-\cdots
\end{align*}

 Applying this to (\ref{condition}), the resultant expression determines
 $n$ differential equations in $m$ functions,   $C'_{i}=F(C_{1},\dots,C_{i-1})$ $(i
\geq 2)$ for an integral curve of the superhorizontal
distribution, where $n=\sum_{i\geq 2}\dim \mathfrak{g}^{I}_i$
 and $m=\sum_{i\geq 1}\dim \mathfrak{g}^{I}_i$.
 These equations can be
 written in local coordinates corresponding to $\bigoplus_{i>0}\mathfrak{g}^{I}_{i}$ as $1$-forms of $m$-variables
giving in this way a Pfaffian system $S$ on $\mathbb{C}^{m}$ with
$\text{rank}(S)=n$.

\begin{prop}
The Pfaffian system $S$ just obtained satisfies the following
properties:
\begin{itemize}
 \item  $\text{class}(S)=\sum_{i\geq 1}\dim \mathfrak{g}^{I}_i=m$.
\item  The first derived system $S'$ of $S$ corresponds to:
\begin{align*}
G^{\mathbb{C}}\times_{P_I}\left((\mathfrak{p}_I\oplus\mathfrak{g}^{I}_{1}\oplus\mathfrak{g}^{I}_{2})/\mathfrak{p}_I\right).
\end{align*}
\end{itemize}
\begin{proof}  Notice first that $\dim \mathcal{A}(\mathcal{I}(SH))= 0 $
is equivalent to $\mathcal{A}(\mathcal{I}(SH))=L(S)\cap
\Gamma(SH)= 0$ , \textit{i.e.}\ if $X\in \Gamma(SH)$ satisfies
\begin{align}
[X,\Gamma(SH)]\subseteq\Gamma(SH), \label{simetria}
\end{align}
then $X=0$.  We show that this is the case. Take
 $X\in SH$ satisfying (\ref{simetria}). Then $X$ is expressed as
 $X=[(g,[Z])]$ with $g\in G^{\mathbb{C}}$ and $[Z]\in \mathfrak{g}_{1}^{I}/
 \mathfrak{p}_{I}$.  Now under the isomorphism $f$, the condition
 (\ref{simetria}) is equivalent to having $Z\in
 \mathfrak{g}_{1}^{I}$ such that $[Z,\mathfrak{g}_{1}^{I}]\in
 \mathfrak{g}_{1}^{I}$. But since $[\mathfrak{g}_{i}^{I},\mathfrak{g}_{j}^{I}]\in
 \mathfrak{g}_{i+j}^{I}$ we must have $Z=0$.\\

 Similarly,
$SH^\prime=\ker(S^\prime)$ is spanned  by:
$$\Gamma(SH)+[\Gamma(SH),\Gamma(SH)],$$ which under the isomorphism
$f$ is equivalent to:
\begin{align*}
SH^\prime&\cong
G^{\mathbb{C}}\times_{P_I}\left(\mathfrak{p}_I\oplus(\mathfrak{g}^I_1+[\mathfrak{g}_1^I,\mathfrak{g}^I_1])/\mathfrak{p}_I\right)\\
&\cong G^{\mathbb{C}}\times_{P_I}\left(\mathfrak{p}_I\oplus\mathfrak{g}^I_1\oplus\mathfrak{g}_2^I/\mathfrak{p}_I\right).\\
\end{align*}
\end{proof}
\end{prop}

  Now we concentrate entirely on the case when $m\leq 5$ and
compute all Pfaffian systems that arises from twistor fibrations
of semisimple Lie algebras by using the above construction and
identifying the obtained system -- after a suitable change of
local coordinates -- with a model in the classification of low
dimensional Pfaffian system mentioned in section 1 (tables 1, 2 and 3).\\

Theorem (\ref{mainteo}) collects these results in a detailed list.
Since the proof is purely computational, in order to illustrate
how the calculation is carried out we present three representative
examples:


\subsection{$\pi:F_3(\mathbb{C}) \rightarrow \mathbb{C}P^2$.}

Let $G=SU_3$, then $G^{\mathbb{C}}=SL_3\mathbb{C}$ with Lie
algebra $\mathfrak{sl}_{3}\mathbb{C}$ consisting of all
endomorphisms of $ \mathbb{C}^{3}$ with zero trace. Let $
\mathfrak{t}^{\mathbb{C}}=\{H\in\mathfrak{sl}_3\mathbb{C}\ |
H=\text{diag}(a_{1},a_{2},a_{3})\}$ be the corresponding Cartan
subalgebra. \\

The set of roots is then given by: $\Delta=\{ \pm(r_i-r_j) \ |
1\leq i,j \leq 3\, \ i\neq j\}$ where $r_{i}$ denotes the linear
functional $r_i:\mathfrak{t}^{\mathbb{C}} \rightarrow \mathbb{C}$
defined by $r_{i}(H)=a_{i}$.   In this case the set of simple
roots is $\Pi=\{\alpha_1=r_{1}-r_{2},\alpha_2=r_{2}-r_{3}\}$ and
the corresponding root spaces are:
 \begin{align}
\mathfrak{g}_{\alpha_1}=\mathbb{C}\cdot
E_{12}~,~\mathfrak{g}_{\alpha_2}=\mathbb{C}\cdot
E_{23}~,~\mathfrak{g}_{\alpha_1+\alpha_2}=\mathbb{C}\cdot
E_{13}\label{rootspacesoff3}
\end{align}
 where $E_{ij}$ stands for the matrix whose
$(i,j)$ entry is 1  and all others 0. \\

If we choose $I=\Pi$, then:
$$\mathfrak{g}_1^I=\left\{\left.\left(
            \begin{array}{@{\,}ccc@{\,}}
                0 & w_1 & 0 \\
                0 & 0 & w_2 \\
                0 & 0 & 0 \\
            \end{array}
        \right) ~\right|~ w_1,w_2 \in \mathbb{C} \right\} ~,~
\mathfrak{g}_2^I=\left\{\left.\left(
            \begin{array}{@{\,}ccc@{\,}}
                0 & 0 & w_3 \\
                0 & 0 & 0 \\
                0 & 0 & 0 \\
            \end{array}
        \right) ~\right|~ w_3 \in \mathbb{C} \right\}.$$
and the corresponding parabolic subalgebra is:
\begin{align*}
\mathfrak{p}_I=\mathfrak{t}^{\mathbb{C}}\oplus\mathfrak{g}_{-\alpha_1}\oplus\mathfrak{g}_{-\alpha_2}\oplus\mathfrak{g}_{-\alpha_3}
          =\left\{\left(
            \begin{array}{@{\,}ccc@{\,}}
                * & 0 & 0 \\
                * & * & 0 \\
                * & * & * \\
            \end{array}
        \right) \in \mathfrak{sl}_3\mathbb{C} \right\}
\end{align*}
with parabolic subgroup:
\begin{align}
P_I=\{g\in
G^{\mathbb{C}}~|~Ad(g)\mathfrak{p}_I\subseteq\mathfrak{p}_I\}
   =\left\{\left(
            \begin{array}{@{\,}ccc@{\,}}
                * & 0 & 0 \\
                * & * & 0 \\
                * & * & * \\
            \end{array}
        \right) \in SL_3\mathbb{C} \right\}\label{parabolicf3}
\end{align}

Since $H=P_I\cap G\cong S(U_1\times U_1\times U_1)$ the
generalized flag manifold $G/H$ is  $SU_3/S(U_1\times U_1\times
U_1)\cong F_3(\mathbb{C})$.\\

The canonical element is calculated as follows: Taking
$\xi\in\mathfrak{t}^\mathbb{C}$ as
$$\xi=\text{diag}(\sqrt{-1}a_1,\sqrt{-1}a_2,\sqrt{-1}a_3)$$ with
$a_i\in\mathbb{R}$, $\sum a_i=0$, since
$\alpha_i(\xi)=\sqrt{-1}n_{I}(\alpha_i)$ we get
 $\alpha_1(\xi)=\alpha_2(\xi)=1$ and $\alpha_3(\xi)=2$, which implies
$a_1=1,a_2=0,a_3=-1$.
 \\

We have $g=\exp\pi\xi=\text{diag}(e^{\sqrt{-1}\pi},1
,e^{-\sqrt{-1}\pi} )$ and therefore the set of fixed points of
$\sigma_g$ is:

\begin{align*}
 G^{\sigma_g}&=\left\{ \left(
            \begin{array}{@{\,}ccc@{\,}}
                a_{11} & 0 & a_{13} \\
                0 & a_{22} & 0 \\
                a_{31} & 0 & a_{33} \\
            \end{array}
        \right) \in SU_3 ~|~ a_{22}\in U_1 ~,~ \left(
                            \begin{array}{@{\,}cc@{\,}}
                                a_{11} & a_{13} \\
                                a_{31} & a_{33} \\
                            \end{array}
                            \right)\in U_2 \right\}\\
                            &\cong S(U_1\times U_2).
 \end{align*}

Here since $H\subset G^{\sigma_g}$, setting  $K=G^{\sigma_g}$ we
can write:
\begin{align*}
G/K\cong SU_3/S(U_1\times U_2)\cong \mathbb{C}P^2.
 \end{align*}

Thus in this case the twistor fibration is precisely
$\pi:F_3(\mathbb{C})\rightarrow \mathbb{C}P^2$.\\

 Now, the integral curve is constructed defining the map $C:\mathbb{C}
\rightarrow \bigoplus_{i=1}^2\mathfrak{g}_i^I$ as:
\begin{align*}
C(z)=\left(
            \begin{array}{@{\,}ccc@{\,}}
                0 & a(z) & c(z) \\
                0 & 0 & b(z) \\
                0 & 0 & 0 \\
            \end{array}
        \right).
\end{align*}

Since
 \begin{align*}
\mathfrak{g}_1^I\oplus \mathfrak{g}_2^I&=\left\{\left.\left(
            \begin{array}{@{\,}ccc@{\,}}
                0 & w_1 & w_3 \\
                0 & 0 & w_2 \\
                0 & 0 & 0 \\
            \end{array}
        \right) ~\right|~ w_1,w_2,w_3 \in \mathbb{C} \right\}
 \end{align*}

we have:
\begin{align*}\exp(-C)(\exp C)^\prime=\left(
            \begin{array}{@{\,}ccc@{\,}}
                0 & a^\prime(z) & c^\prime(z)-\frac{1}{2}a(z)b^\prime(z)+\frac{1}{2}a^\prime(z)b(z) \\
                0 & 0 & b^\prime(z) \\
                0 & 0 & 0 \\
            \end{array}
        \right)
 \end{align*}

After applying the condition (\ref{condition}) we obtain the
following ordinary differential equation:
\begin{align*}
c^\prime(z)-\frac{1}{2}a(z)b^\prime(z)+\frac{1}{2}a^\prime(z)b(z)=0.
\end{align*}
In terms of local coordinates of
$\exp\bigoplus_{i=0}^2\mathfrak{g}_i^I$ this becomes
\begin{align*}
dw_3-\frac{1}{2}w_1dw_2+\frac{1}{2}w_2dw_1=0.
\end{align*}
Since
$$dw_3-\frac{1}{2}w_1dw_2+\frac{1}{2}w_2dw_1=d(w_3-\frac{1}{2}w_1w_2)+w_2dw_1,$$
the  change of coordinates
$$(z_1,z_2,z_3)=(w_3-\frac{1}{2}w_1w_2,w_2,w_1)$$
identifies this Pfaffian system with:
\begin{align*}
 S^2_3=\{dz_1+z_2dz_3\}.
\end{align*}

\subsection{$\pi:SO(5)/U(1)\times SO(2) \rightarrow SO(5)/SO(2)\times SO(3)$.}

~~ Let  $G=SO(5)$. By definition the corresponding Lie algebra is
given by:
\begin{align*}
\mathfrak{so}_{5}^\mathbb{C}=\{ X \in
Hom(\mathbb{C}^5,\mathbb{C}^5) ~|~ (Xv,w)+(v,Xw)=0 \}.
\end{align*}
 After choosing a basis:
\begin{align*}
v_1&=e_1+\sqrt{-1}e_4, ~~~ v_2=e_2+\sqrt{-1}e_5,~~~v_3=e_3,\\
v_4&=e_2-\sqrt{-1}e_5,~~~ v_5=e_1-\sqrt{-1}e_4.
\end{align*}
where $\{e_i\}_{i=1}^{5}$ denotes the canonical basis of
$\mathbb{C}^{5}$, a suitable matrix representation is obtained:

\begin{align*}
\mathfrak{so}(5,\mathbb{C})=\left\{\left.\left(
            \begin{array}{@{\,}ccccc@{\,}}
                r_1 & x & t_1 & y & 0 \\
                \tilde{x} & r_2 & t_2 & 0 & -y \\
                2\tilde{t_1} & 2\tilde{t_2} & 0 & -2t_2 & -2t_1
                \\
                \tilde{y} & 0 & -\tilde{t_2} & -r_2 & -x \\
                0 & -\tilde{y} & -\tilde{t_1} & -\tilde{x} & -r_1
            \end{array}
        \right) ~\right|~
        r_1,r_2,x,y,t_1,t_2,\tilde{x},\tilde{y},\tilde{t_1},\tilde{t_2}\in\mathbb{C}\right\}.
\end{align*}

In this case the set of roots is: $\Delta=\{\pm(r_i\pm r_j),\pm
r_k ~|~1\leq i,j,k\leq 2 , i\neq j\}$ with simple roots
$\Pi=\{\alpha_1=r_1-r_2,\alpha_2=r_2\}$, and the root spaces are:
\begin{align*}
    \mathfrak{g}_{\alpha_1}&=\mathbb{C} \cdot (E_{12} -E_{45}), \ \
     \mathfrak{g}_{\alpha_2}=\mathbb{C} \cdot (E_{23} -2E_{34}),
     \\
    \mathfrak{g}_{\alpha_1+\alpha_2}&=\mathbb{C} \cdot (E_{13}
    -2E_{35}), \ \
    \mathfrak{g}_{\alpha_1+2\alpha_2}=\mathbb{C} \cdot (E_{14}
    -E_{25}).
\end{align*}

 Setting $I=\Pi$, we have:
 \begin{align*}
\mathfrak{g}_{1}^I=
\mathfrak{g}_{\alpha_1}\oplus\mathfrak{g}_{\alpha_2} , ~~~
\mathfrak{g}_2^I=\mathfrak{g}_{\alpha_1+\alpha_2} ,~~~
\mathfrak{g}_3^I=\mathfrak{g}_{\alpha_1+2\alpha_2}.
\end{align*}

The map $C:\mathbb{C} \rightarrow \bigoplus_{i>0}\mathfrak{g}_i^I
$ can be defined as:

\begin{align*}
C(z)=\left(
            \begin{array}{@{\,}ccccc@{\,}}
                0 & x(z) & t_1(z) & y(z) & 0 \\
                0 & 0 & t_2(z) & 0 & -y(z) \\
                0 & 0 & 0 & -2t_2(z) & -2t_1(z) \\
                0 & 0 & 0 & 0 & -x(z) \\
                0 & 0 & 0 & 0 & 0
            \end{array}
        \right)
\end{align*}

In the same way as above, calculating $(\exp -C)(\exp C)^\prime$
under the condition (\ref{condition}), we obtain the following
system of differential equations:
\begin{align*}
t_1^\prime-\frac{1}{2}xt_2^\prime+\frac{1}{2}t_2x^\prime=0,~~~~
y^\prime+t_1t_2^\prime-t_2t_1^\prime+\frac{1}{3}t_2xt_2^\prime-\frac{1}{3}t_2^2x^\prime=0.
\end{align*}

In terms of local coordinates these are:
$dt_1-\frac{1}{2}xdt_2+\frac{1}{2}t_2dx=0$ and
$dy+t_1dt_2-t_2dt_1+\frac{1}{3}t_2xdt_2-\frac{1}{3}t_2^2dx=0$. The
associated Pfaffian system $S=\{\omega_1,\omega_2\}$ is given by:
\begin{align*}
 \omega_1:&=dt_1-\frac{1}{2}xdt_2+\frac{1}{2}t_2dx,\\
\omega_2:&=dy+t_1dt_2-t_2dt_1+\frac{1}{3}t_2xdt_2-\frac{1}{3}t_2^2dx.
\end{align*}

Notice that $S$ is not completely integrable and that the derived
system is $S^\prime=\{\omega_2\}$ since
$d\omega_2\wedge\omega_1\wedge\omega_2=0$. Moreover, since
$d\omega_2\wedge\omega_2\neq0$ we have
$\textrm{class}(S^\prime)=3$ and the system must be equivalent to
$S^7_4$. We shall find an explicit change of coordinates. First,
 $\omega_{2}$ can be expressed as:
$$\omega_2=d(y-t_1t_2-\frac{1}{3}xt_2^2)+(2t_1+xt_2)dt_2,$$ to obtain:
 \begin{align*}
(z_1,z_2,z_3)=(y-t_1t_2-\frac{1}{3}xt_2^2,~~2t_1+xt_2,~~t_2).
 \end{align*}

On the other hand:
$$\omega_1=\frac{1}{2}(d(2t_1+xt_2)+(-2x)dt_2)$$
hence the change of coordinates:
$$(z_1,z_2,z_3,z_4)=(y-t_1t_2-\frac{1}{3}xt_2^2,~~2t_1+xt_2,~~t_2,~~-2x)$$
shows that the Pfaffian system is :
\begin{align*}
S^7_4=\{dz_1+z_2dz_3,~~ dz_2+z_4dz_3\}.
\end{align*}

\subsection{$\pi:G_2/U(2) \rightarrow G_2/SO(4)$.}

    ~ In this example we consider $G_2=\{X\in SO(7,\mathbb{R}) ~|~X(v\times w)=Xv\times Xw
  \}$, with Lie algebra:
\begin{align*}
\mathfrak{g}_2^{\mathbb{C}}=\{ X \in \mathfrak{so}_7^\mathbb{C}
~|~X(v\times w) =Xv\times w +v\times Xw \}.
\end{align*}

As for the last example, in order to obtain a matrix
representation for  $SO(7)$ we choose a basis:
\begin{align*}
v_1&=e_1+\sqrt{-1}e_5,~~v_2=e_2+\sqrt{-1}e_6,~~v_3=e_3+\sqrt{-1}e_7,~~v_4=e_4,\\
v_5&=e_3-\sqrt{-1}e_7,~~v_6=e_2-\sqrt{-1}e_6,~~v_7=e_1-\sqrt{-1}.
\end{align*}

Then under this representation any element of
$\mathfrak{g}_2^{\mathbb{C}}$ can be expressed as:

$$\left(
            \begin{array}{@{\,}ccccccc@{\,}}
                r_2+r_3 & x_1 & x_2 & y_3 & y_2 & y_1 & 0 \\
                \tilde{x_1} & r_2 & x_3 & -x_2 & y_3 & 0 & -y_1 \\
                \tilde{x_2} & \tilde{x_3} & r_3 & x_1 & 0 & -y_3 & -y_2 \\
                2\tilde{y_2} & -2\tilde{x_2} & 2\tilde{x_1} & 0 & -2x_1 & 2x_2 & -2y_3 \\
                \tilde{y_2} & \tilde{y_3} & 0 & -\tilde{x_1} & -r_3 & -x_3 & -x_2  \\
                \tilde{y_1} & 0 & -\tilde{y_3} & \tilde{x_2} & -\tilde{x_3} & -r_2 & -x_1 \\
                0 & -\tilde{y_1} & -\tilde{y_2} & -\tilde{y_3} & -\tilde{x_2} & -\tilde{x_1} & -r_2-r_3
            \end{array}
        \right)$$

 In this case the set of roots is:
 \begin{align*}
\Delta=\{\pm\alpha_1,\pm\alpha_2,\pm(\alpha_1+\alpha_2),\pm(2\alpha_1+\alpha_2),\pm(3\alpha+\alpha_2),\pm(3\alpha_1+2\alpha_2)\}
\end{align*}
with simple roots  $\Pi=\{\alpha_1,\alpha_2\}$, where
$\alpha_1=r_3$ and $\alpha_2=r_2-r_3$. Here the root spaces are:
\begin{align*}
    \mathfrak{g}_{\alpha_1}&=\mathbb{C} \cdot (E_{12} +E_{34}-2E_{45}-E_{67}),\\
    \mathfrak{g}_{\alpha_2}&=\mathbb{C} \cdot (E_{23} -E_{56}),\\
    \mathfrak{g}_{\alpha_1+\alpha_2}&=\mathbb{C} \cdot (E_{13}
    -E_{24}+2E_{46}-E_{57}),\\
            \mathfrak{g}_{2\alpha_1+\alpha_2}&=\mathbb{C} \cdot (E_{14} +E_{25}-E_{36}-2E_{47}),\\
    \mathfrak{g}_{3\alpha_1+\alpha_2}&=\mathbb{C} \cdot
    (E_{15}-E_{37}),\\
            \mathfrak{g}_{3\alpha_1+2\alpha_2}&=\mathbb{C} \cdot
            (E_{16}-E_{27}).
\end{align*}

Taking  $I=\{\alpha_1\}$ we have :
\begin{align*}
\mathfrak{g}_{1}^I=
\mathfrak{g}_{\alpha_1}\oplus\mathfrak{g}_{\alpha_1+\alpha_2} ,~~
\mathfrak{g}_2^I=\mathfrak{g}_{2\alpha_1+\alpha_2} ,~~
\mathfrak{g}_3^I=\mathfrak{g}_{3\alpha_1+\alpha_2}\oplus\mathfrak{g}_{3\alpha_1+2\alpha_2}.
\end{align*}

Thus the map $C:\mathbb{C} \rightarrow
\bigoplus_{i>0}\mathfrak{g}_i^I$ takes the form:

\begin{align*}
C(z)=\left(
            \begin{array}{@{\,}ccccccc@{\,}}
                0 & x_1(z) & x_2(z) & y_3(z) & y_2(z) & y_1(z) & 0 \\
                0&0 & 0 & -x_2(z) & y_3(z) & 0 & -y_1(z) \\
                0&0&0 & x_1(z) & 0 & -y_3(z) & -y_2(z) \\
                0&0&0 & 0 & -2x_1(z) & 2x_2(z) & -2y_3(z) \\
                0&0&0&0&0 & 0 & -x_2(z)  \\
                0&0&0&0&0&0 & -x_1(z) \\
                0 & 0&0&0&0&0&0
            \end{array}
        \right)
\end{align*}

This time the condition of the integral curve to the superhorizontal distribution gives the following system of differential
equations:
\begin{align*}
y_3^\prime-x_2x_1^\prime+x_1x_2^\prime=&0,\\
y_2^\prime-\frac{3}{2}x_1y_3^\prime+\frac{3}{2}y_3x_1^\prime-x_1^2x_2^\prime+x_1x_2x_1^\prime=&0,\\
y_1^\prime+\frac{3}{2}x_2y_3^\prime-\frac{3}{2}y_3x_2^\prime+x_1x_2x_2^\prime-x_2^2x_1^\prime=&0.
\end{align*}

After expressing them in terms of local coordinates we set

\begin{align*}
\omega_1:&=dy_3-x_2dx_1+x_1dx_2,\\
\omega_2:&=dy_2-\frac{3}{2}x_1dy_3+\frac{3}{2}y_3dx_1-x_1^2dx_2+x_1x_2dx_1,\\
\omega_3:&=dy_1+\frac{3}{2}x_2dy_3-\frac{3}{2}y_3dx_2+x_1x_2dx_2-x_2^2dx_1.
\end{align*}
Then the induced Pfaffian system is $S=\{\omega_1,
~\omega_2,~\omega_3\}$ which is not completely integrable.  Now,
for $i=2,3$  we have that:
\begin{align*}
d\omega_i\wedge\omega_1\wedge\omega_2\wedge\omega_3=0.
\end{align*}
 Therefore, the derived system is $S^\prime=\{\omega_2,\omega_3\}$.
Moreover  since $d\omega_2=-3\omega_1\wedge dx_2$ and
$d\omega_3=3\omega_1\wedge dx_1$ we have
$\textrm{class}(S^\prime)=5$, and then the system
 $S$ is identified to be $S^{11}_5(f)$.\\

Now $\omega_3$ can be written as:
\begin{align*}
\omega_3&=d(y_2-\frac{3}{2}x_1y_3-x_1^2x_2)+(3y_3+3x_1x_2)dx_1
\end{align*}
hence:
\begin{align*}
(w_1,w_3,w_4)=(y_2-\frac{3}{2}x_1y_3-x_1^2x_2,~~3y_3+3x_1x_2,~~x_1)
\end{align*}

In the same way:
\begin{align*}
\omega_2&=d(y_1-\frac{3}{2}x_2y_3-x_1x_2^2)+(x_2)d(3y_3+3x_1x_2)-3x_2^2dx_1,\\
\omega_1&=\frac{1}{3}(d(3y_3+3x_1x_2)-6x_2dx_1).
\end{align*}
implies:
\begin{align*}
(w_1,w_2,w_3,&w_4,w_5)=\\
&(y_2-\frac{3}{2}x_1y_3-x_1^2x_2,~~y_1-\frac{3}{2}x_2y_3-x_1x_2^2,~~3y_3+3x_1x_2,~~x_1,~~x_2).
\end{align*}

Therefore we must have $S^{11}_5(f)$ with $f=-3w_5^2$. But in
general, the  change of coordinates corresponding to
$S^{11}_5(c\cdot w_5^2)$ is:
\begin{align*}
(z_1,z_2,z_3,z_4,z_5)=(w_1,w_2,\frac{1}{\sqrt[3]{c}}w_3,\sqrt[3]{c}w_4,\sqrt[3]{c}w_5).
\end{align*}

Thus we can identify our Pfaffian system  to be of the form
$S^{11}_5(z_5^2)$. Explicitly this can be expressed as
\begin{align*}
S^{11}_5(z_5^2) = \{ dz_1+z_3dz_4 ,~~ dz_2+z_5dz_3+z_5^2dz_4 ,~~
dz_3+2z_5dz_4 \}.
\end{align*}
under the change of coordinates:

\begin{align*}
z_1&=y_2-\frac{3}{2}x_1y_3-x_1^2x_2,\\
z_2&=y_1-\frac{3}{2}x_2y_3-x_1x_2^2,\\
z_3&=-\frac{3}{\sqrt[3]{3}}y_3-\frac{3}{\sqrt[3]{3}}x_1x_2,\\
z_4&=-\sqrt[3]{3}x_1,\\
z_5&=-\sqrt[3]{3}x_2.
\end{align*}

After the examples finally we are in position to state the main
theorem of this section:

\begin{theo} \label{mainteo}
The Pfaffian systems of at most five variables that arise from superhorizontal distributions of
twistor fibrations of semisimple Lie algebras are given explicitly
in table 4.
\begin{proof} By direct calculation.\\\end{proof}
\end{theo}

\begin{itemize}
\item [\textbf{Table 4.}]  Low dimensional Pfaffian systems that
arise from twistor fibrations.
\end{itemize}
    \rotatebox{270}{
    \begin{tabular}{|c|c|c|}
\hline
        Lie group &
        decomposition of $T_{1,0}G/H$ &
        Pfaffian system     \\ \hline &&\\[-3pt]
        $SU_3$  &    &  \\[-18pt]
        &&$S^2_3=\{dz_1+z_2dz_3\}$\\
        &&{\scriptsize$(w_3-\frac{1}{2}w_1w_2,w_2,w_1)$}\\[-25pt]
    &$\overbrace{\mathfrak{g}_{\alpha}\oplus\mathfrak{g}_{\beta}}^{\mathfrak{g}^I_1}\bigoplus\overbrace{\mathfrak{g}_{\alpha+\beta}}^{\mathfrak{g}_2^I}$          & \\ \hline &&\\[28pt]
        $SU_4$ & &\\[-48pt]
        &&$S^{15}_5(0)=\{ dz_1+z_3dz_4 , dz_2+z_5dz_3\}$\\
        &&{\scriptsize$(w_5+\frac{1}{2}w_1w_3,w_4-\frac{1}{2}w_1w_2,w_1,-w_3,w_2)$}\\[-26pt]
    &   $\overbrace{\mathfrak{g}_{\alpha}\oplus\mathfrak{g}_{\beta}\oplus\mathfrak{g}_{\beta+\gamma}}^{\mathfrak{g}^I_1}\bigoplus\overbrace{\mathfrak{g}_{\alpha+\beta}\oplus\mathfrak{g}_{\alpha+\beta+\gamma}}^{\mathfrak{g}_2^I}$ &   \\ \cline{2-3} &&\\[-7pt]
        &&$S^{15}_5(0)=\{ dz_1+z_3dz_4 , dz_2+z_5dz_3\}$\\
        &&{\scriptsize$(w_5-\frac{1}{2}w_2w_3,w_4+\frac{1}{2}w_1w_2,w_2,w_3,-w_1)$}\\[-26pt]
    &   $\overbrace{\mathfrak{g}_{\beta}\oplus\mathfrak{g}_{\gamma}\oplus\mathfrak{g}_{\alpha+\beta}}^{\mathfrak{g}^I_1}\bigoplus\overbrace{\mathfrak{g}_{\beta+\gamma}\oplus\mathfrak{g}_{\alpha+\beta+\gamma}}^{\mathfrak{g}_2^I}$ &   \\ \cline{2-3} &&\\[-7pt]
        &&$S^3_5=\{dz_1+z_2dz_3+z_4dz_5\}$\\
        &&{\scriptsize$(w_5-\frac{1}{2}w_1w_4-\frac{1}{2}w_2w_3,w_4,w_1,w_2,w_3)$}\\[-26pt]
            &   $\overbrace{\mathfrak{g}_{\alpha}\oplus\mathfrak{g}_{\gamma}\oplus\mathfrak{g}_{\alpha+\beta}\oplus\mathfrak{g}_{\beta+\gamma}}^{\mathfrak{g}^I_1}\bigoplus\overbrace{\mathfrak{g}_{\alpha+\beta+\gamma}}^{\mathfrak{g}_2^I}$        &  \\ \hline &&\\[14pt]
        $SO_5$ &&\\[-34pt]
        &&$S^2_3=\{dz_1+z_2dz_3\}$\\
        &&{\scriptsize$(w_3-w_1w_2,2w_2,w_1)$}\\[-26pt]
    &   $\overbrace{\mathfrak{g}_{\beta}\oplus\mathfrak{g}_{\alpha+\beta}}^{\mathfrak{g}^I_1}\bigoplus\overbrace{\mathfrak{g}_{\alpha+2\beta}}^{\mathfrak{g}_2^I}$       &      \\ \cline{2-3} &&\\[-7pt]
        &&$S^7_4=\{ dz_1+z_2dz_3 , dz_2+z_4dz_3 \}$\\
        &&{\scriptsize$(w_4-w_2w_3-\frac{1}{3}w_1w_2^2,2w_3+w_1w_2,w_2,-2w_1)$}\\[-26pt]
            &   $\overbrace{\mathfrak{g}_{\alpha}\oplus\mathfrak{g}_{\beta}}^{\mathfrak{g}^I_1}\bigoplus\overbrace{\mathfrak{g}_{\alpha+\beta}}^{\mathfrak{g}_2^I}\bigoplus\overbrace{\mathfrak{g}_{\alpha+2\beta}}^{\mathfrak{g}_3^I}$ &    \\ \hline &&\\[28pt]
        $SO_6$ &&\\[-48pt]
        &&$S^3_5=\{dz_1+z_2dz_3+z_4dz_5\}$\\
        &&{\scriptsize$(w_5-\frac{1}{2}w_2w_3-\frac{1}{2}w_1w_4,w_3,w_2,w_4,w_1)$}\\[-26pt]
    &   $\overbrace{\mathfrak{g}_{\beta}\oplus\mathfrak{g}_{\gamma}\oplus\mathfrak{g}_{\alpha+\beta}\oplus\mathfrak{g}_{\alpha+\gamma}}^{\mathfrak{g}^I_1}\bigoplus\overbrace{\mathfrak{g}_{\alpha+\beta+\gamma}}^{\mathfrak{g}_2^I}$ &  \\ \cline{2-3} &&\\[-7pt]
        &&$S^{15}_5(0)=\{ dz_1+z_3dz_4 , dz_2+z_5dz_3\}$\\
        &&{\scriptsize$(w_4+\frac{1}{2}w_1w_2,w_5-\frac{1}{2}w_2w_3,-w_1,w_2,w_3)$}\\[-26pt]
            &   $\overbrace{\mathfrak{g}_{\alpha}\oplus\mathfrak{g}_{\gamma}\oplus\mathfrak{g}_{\alpha+\beta}}^{\mathfrak{g}^I_1}\bigoplus\overbrace{\mathfrak{g}_{\alpha+\gamma}\oplus\mathfrak{g}_{\alpha+\beta+\gamma}}^{\mathfrak{g}_2^I}$ &        \\ \cline{2-3} &&\\[-7pt]
        &&$S^{15}_5(0)=\{ dz_1+z_3dz_4 , dz_2+z_5dz_3\}$\\
        &&{\scriptsize$(w_4+\frac{1}{2}w_1w_2,w_5-\frac{1}{2}w_2w_3,-w_1,w_2,w_3)$}\\[-26pt]
            &   $\overbrace{\mathfrak{g}_{\alpha}\oplus\mathfrak{g}_{\beta}\oplus\mathfrak{g}_{\alpha+\gamma}}^{\mathfrak{g}^I_1}\bigoplus\overbrace{\mathfrak{g}_{\alpha+\beta}\oplus\mathfrak{g}_{\alpha+\beta+\gamma}}^{\mathfrak{g}_2^I}$ &     \\ \hline &&\\[14pt]
        $Sp_2$ &&\\[-34pt]
        &&$S^2_3=\{dz_1+z_2dz_3\}$\\
        &&{\scriptsize$(w_3-w_1w_2,2w_2,w_1)$}\\[-26pt]
    &   $\overbrace{\mathfrak{g}_{\alpha}\oplus\mathfrak{g}_{\alpha+\beta}}^{\mathfrak{g}^I_1}\bigoplus\overbrace{\mathfrak{g}_{2\alpha+\beta}}^{\mathfrak{g}_2^I}$          &      \\ \cline{2-3} &&\\[-7pt]
        &&$S^7_4=\{ dz_1+z_2dz_3 , dz_2+z_4dz_3 \}$\\
        &&{\scriptsize$(w_4-w_1w_3+\frac{1}{3}w_1^2w_2,2w_3-w_1w_2,w_1,2w_2)$}\\[-26pt]
            &   $\overbrace{\mathfrak{g}_{\alpha}\oplus\mathfrak{g}_{\beta}}^{\mathfrak{g}^I_1}\bigoplus\overbrace{\mathfrak{g}_{\alpha+\beta}}^{\mathfrak{g}_2^I}\bigoplus\overbrace{\mathfrak{g}_{\alpha+2\beta}}^{\mathfrak{g}_3^I}$ &    \\ \hline &&\\[-3pt]
        $Sp_3$ &     &   \\[-18pt]
        &&$S^3_5=\{dz_1+z_2dz_3+z_4dz_5\}$\\
        &&{\scriptsize$(w_5-w_1w_4-w_2w_3,2w_4,w_1,2w_3,w_2)$}\\[-26pt]
        &$\overbrace{\mathfrak{g}_{\alpha}\oplus\mathfrak{g}_{\alpha+\beta}\oplus\mathfrak{g}_{\alpha+\beta+\gamma}\oplus\mathfrak{g}_{\alpha+2\beta+\gamma}}^{\mathfrak{g}^I_1}\bigoplus\overbrace{\mathfrak{g}_{2\alpha+2\beta+\gamma}}^{\mathfrak{g}_2^I}$&\\ \hline
&&\\[14pt]
        $G_2$ &&\\[-34pt]
        &&$S^{11}_5(z_5^2)= \{ dz_1+z_3dz_4 , dz_2+z_5dz_3+z_5^2dz_4 , dz_3+2z_5dz_4 \}$\\      &&{\scriptsize$(w_4-\frac{3}{2}w_1w_3-w_1^2w_2,w_5-\frac{3}{2}w_2w_3-w_1w_2^2,-\frac{3}{\sqrt[3]{3}}w_3-\frac{3}{\sqrt[3]{3}}w_1w_2,-\sqrt[3]{3}w_1,-\sqrt[3]{3}w_2)$} \\[-27pt]
    &   $\overbrace{\mathfrak{g}_{\alpha}\oplus\mathfrak{g}_{\alpha+\beta}}^{\mathfrak{g}^I_1}\bigoplus\overbrace{\mathfrak{g}_{2\alpha+\beta}}^{\mathfrak{g}_2^I}\bigoplus\overbrace{\mathfrak{g}_{3\alpha+\beta}\oplus\mathfrak{g}_{3\alpha+2\beta}}^{\mathfrak{g}_3^I}$       &      \\ \cline{2-3} &&\\[-7pt]
        &&$S^3_5=\{dz_1+z_2dz_3+z_4dz_5\}$\\
        &&{\scriptsize$(w_5-\frac{3}{2}w_2w_3-\frac{1}{2}w_1w_4,3w_2,w_3,w_4,w_1)$}\\[-26pt]
            &   $\overbrace{\mathfrak{g}_{\beta}\oplus\mathfrak{g}_{\alpha+\beta}\oplus\mathfrak{g}_{2\alpha+\beta}\oplus\mathfrak{g}_{3\alpha+\beta}}^{\mathfrak{g}^I_1}\bigoplus\overbrace{\mathfrak{g}_{3\alpha+2\beta}}^{\mathfrak{g}_2^I}$      &              \\ \hline
\end{tabular}}

\newpage

\textbf{Remark:} Since  $\mathfrak{so}_{5}$ is isomorphic to
$\mathfrak{sp}_{2}$ and $\mathfrak{so}_{6}$  to
$\mathfrak{su}_{4}$, the corresponding  Pfaffian systems must be
 equivalent. However since in each case the details of the
calculation are different and for the sake of confirmation their
description is also included in the above list.\\

\section{Application: infinitesimal symmetries.}

 The remainder of this note presents an example of how the
twistor Pfaffian systems discussed above are related to the
corresponding Lie algebra of infinitesimal
symmetries, which in general is infinite dimensional.\\

It is well known that $\mathfrak{g}^\mathbb{C}\subseteq L(S)$,
because the superhorizontal distribution is a
$G^\mathbb{C}-$invariant distribution, \textit{i.e.}\
$G^\mathbb{C}$ acts on $G^\mathbb{C}/P_{I}$ by left translations
preserving $SH$. More precisely,  for any
$X\in\mathfrak{g}^{\mathbb{C}}$ the infinitesimal symmetries of
$SH$ are given  in a canonical way by the vector fields:
 \begin{align}X^{*}_{[g]}= \frac{d}{dt}[\exp(tX)g ]|_{t=0},
\label{vectorfield}
 \end{align}
with $[g]\in B$, the big cell of $G^{\mathbb{C}}/P_{I}$.\\

Secondly, E. Cartan observed (see [8]) that the Pfaffian
system $S^{11}_{5}(z^{2}_{5})$ on $\mathbb{C}^{5}$ has the property $L(S^{11}_{5}(z^{2}_{5}))\cong
\mathfrak{g}^{\mathbb{C}}_{2}$ and therefore $\dim L(S^{11}_{5}(z^{2}_{5}))_{\mathbb{C}}=14$.\\

More generally, K. Yamaguchi in [10] considered examples of
regular differential systems, which turn out to agree with the
Pfaffian systems arising from twistor fibrations. For any such
system $S$, the main result of [10] asserts that $L(S)\cong
\mathfrak{g}^{\mathbb{C}}$ except for the following three cases:

\begin{itemize}
\item[(1)]$\mathfrak{g}^{\mathbb{C}}=\mathfrak{g}_{-1}^I\oplus\mathfrak{g}_{0}^I\oplus\mathfrak{g}_{1}^I$.\\
\item[(2)] $\mathfrak{g}^{\mathbb{C}}=\bigoplus_{i=-2}^{2}\mathfrak{g}_i^I$ ~~(If $\dim \mathfrak{g}_{-2}^I=\dim \mathfrak{g}_2^I=1$).\\
\item[(3)] $\mathfrak{g}^{\mathbb{C}}$ is a Lie algebra of type
$A_{l}$ such that $I=\{\alpha_1,\alpha_m\}$, or type $C_{l}$ such
that $ I=\{\alpha_1,\alpha_l\}$.
$(1\hspace{-.1em}<\hspace{-.1em}m\hspace{-.1em}<\hspace{-.1em}l)
$.
\end{itemize}

Notice that all the examples presented in table 4  belong to one of these categories with the exception of  $S^{11}_{5}(z^{2}_{5})$, the Cartan case. \\

  The fact that the Pfaffian system $S^{2}_{3}$ on
$\mathbb{C}^{3}$ originates from \emph{two different} twistor
fibrations (see table 4) give us a  natural explanation of the fact that $\dim
L(S^{2}_{3})=\infty$. In fact, as we have seen in 3.1, the twistor fibration
$$\pi:SU_3/S(U_1\times U_1\times U_1)\cong F_3(\mathbb{C})
\rightarrow SU_3/S(U_1\times U_2)\cong \mathbb{C}P^2,$$  with big
cell $\mathbb{C}^{3}$ has a superhorizontal distribution generated
by $\{\partial_{1},\ \partial_{2}+w_1\partial_{3}\}$ and therefore
a Pfaffian system $S=\{dw_3-w_1dw_2\}$, which is equivalent to
$S_{3}^{2}$ under the change of coordinates:
$(z_{1},z_{2},z_{3})=(w_3,-w_1,w_2).$  By computing
(\ref{vectorfield}) for a basis of the lie algebra
  $\mathfrak{sl}_{3}\mathbb{C}$, a set  of
  8 vector fields is obtained. The Lie algebra $\mathcal{L}'$ spanned by this
  set is a subalgebra of the infinitesimal symmetries of
  $S_{3}^{2}$. Similarly, for the  fibration $$\pi:Sp_{2}/Sp_{1}\times U_{1}
\cong \mathbb{C}P^{3}\rightarrow Sp_{2}/Sp_{1}\times Sp_{1}\cong
S^{4},$$ we have
$SH=\langle\{\partial_{1}-w_2\partial_{3},\partial_{2}+w_1\partial_{3}\}\rangle$
and $S=\{ dw_3+w_2dw_1-w_1dw_2\}$ which under the change of
coordinates $(z_1,z_2,z_3)=(w_3-w_1w_2,2w_2,w_1)$ also corresponds
to $S^{2}_{3}$. In the same way as above,  we can construct a Lie
subalgebra $\mathcal{L}''$ of $L(S_{3}^{2})$ spanned
by  10 vector fields by means of (\ref{vectorfield}).\\

By direct calculation we find that the Lie algebra
$\langle\mathcal{L}',\mathcal{L}''\rangle$ generated by
$\mathcal{L}'$ and $\mathcal{L}''$ is an \textit{infinite
dimensional} Lie subalgebra of $L(S_{3}^{2})$.  In a future paper
we shall discuss further ramifications of this
observation.\\

\textbf{Acknowledgement.} We are grateful to Professors J\"{u}rgen
Berndt and John Bolton for their interest in this work and their
helpful comments, and
to Professor Martin Guest for his invaluable guidance and patience.\\

\begin{itemize}
\item []\textit{Department of Mathematics and Information Sciences.\\
Tokyo Metropolitan University\\
Minami-Ohsawa 1-1, Hachioji-shi,\\
Tokyo 192-0397, Japan.\\\textit{e-mail:}
carrillo-ramirocatalan@c.metro-u.ac.jp.}
\end{itemize}


\begin{thebibliography}{100}

\bibitem [1]{} Awane A., Goze M. \, {\em Pfaffian systems, k-symplectic systems. } Kluwer Academic Publishers, Dordrecht,
 ISBN 0-7923-6373-6, 2000.

\bibitem [2]{} Bryant R.L. \, {\em Lie groups and twistor spaces.} \ Duke Math. J. \textbf{52} 1985, 223-261.


\bibitem [3]{} Bryant R.L.,  Chern S.S.,  Gardner R.B., Goldschmidt H.L., Griffiths P.A., \, {\em Exterior differential systems.} Mathematical Sciences Research Institute Publications 18,
Springer-Verlag, 1991.


\bibitem [4]{}  Burstall F.E., Guest M.A.  \, {\em Harmonic two-spheres in compact symmetric spaces, revisited.} \  Mathematische  Annalen \textbf{309} 1997, 541-572.

\bibitem [5]{}  Burstall F.E., Rawnsley J.H.  \, {\em Twistor theory for Riemmanian symmetric spaces.} \ Lecture Notes in Math. \textbf{1424}. Berlin, Heidelberg: Springer 1990.

\bibitem [6]{} Ca\~nadas M.A., Ruiz C.  \, {\em Pfaffian systems with derived length one. The class of flag systems.} \ Transactions of the American Mathematical Society, \textbf{353}, no.5,
2001,  1755-1766.


\bibitem [7]{  } Humphreys J.E. \, {\em Introduction to Lie algebras and representation theory.} Graduate Text in Mathematics 9. Springer-Verlag.  1972.

\bibitem [8]{} Kumpera A.  \, {\em On the Lie and Cartan theory of invariant differential systems.} \ J. Math. Sci. Univ. Tokyo  \textbf{6},1999, 229-314.

\bibitem [9]{}  Morita S. \, {\em Geometry of differential forms.} \ Iwanami Series of Modern Mathematics AMS 201,  1998.

\bibitem [10]{}  Yamaguchi K. \, {\em Differential systems associated with simple graded Lie algebras.} \  Advanced Studies in Pure Mathematics. Progress in Differential Geometry. 22, 1993,
413-494.\\


\end{thebibliography}
\end{document}